\documentclass[a4,12pt]{amsart}
\usepackage[dvips]{graphicx}
\usepackage{amsmath}
\usepackage{amssymb}
\usepackage{subfigure}
\usepackage{stmaryrd}

\theoremstyle{plain}
\newtheorem{theorem}{Theorem}[section]
\newtheorem{lemma}[theorem]{Lemma}
\newtheorem{corollary}[theorem]{Corollary}
\newtheorem{proposition}[theorem]{Proposition}

\theoremstyle{definition}
\newtheorem{definition}{Definition}
\newtheorem{example}{Example}

\theoremstyle{remark}
\newtheorem{remark}{Remark}
\newtheorem*{notation}{Notation}

\providecommand{\real}{\mathbb{R}}
\providecommand{\rn}{\mathbb{R}^{n}}
\providecommand{\p}{\mathbb{P}}
\providecommand{\pstar}{\mathbb{P}^{\ast}}
\providecommand{\pn}{\mathbb{P}^{n}}
\providecommand{\nat}{\mathbb{N}}

\providecommand{\introp}[1]{\mathop{\rm int}\nolimits #1}
\providecommand{\comp}[1]{{\rm Cmp}\left(#1\right)}
\providecommand{\cocomp}[1]{{\rm Cocmp}\left(#1\right)}
\providecommand{\relhull}[2]{\left[ #1 \right]_{#2}}
\providecommand{\hull}[1]{\left[ #1 \right]}
\providecommand{\pair}[2]{\left\langle #1, #2 \right\rangle}
\providecommand{\sat}[1]{\left\llbracket #1 \right\rrbracket}
\providecommand{\set}[2]{\left\{#1\, |\, #2\right\}}

\begin{document}

\title[Multi-convex sets in real projective spaces]{Multi-convex sets in real projective spaces and their duality}
\author{Takahisa Toda}
\address{Graduate School of Human and Environmental Studies, Kyoto University\\
  Yoshida-nihonmatsu-cho, Sakyo-ku, Kyoto 606-8501, JAPAN}
\email{toda.takahisa@hw3.ecs.kyoto-u.ac.jp}
\date{\today}

\begin{abstract}
We study intersections of projective convex sets in the sense of Steinitz.
In a projective space, an intersection of a nonempty family of convex sets splits into multiple connected components each of which is a convex set.
Hence, such an intersection is called a multi-convex set.
We derive a duality for saturated multi-convex sets:
there exists an order anti-isomorphism between nonempty saturated multi-convex sets in a real projective space and those in the dual projective space.
In discrete geometry and computational geometry, these results allow to transform a given problem into a dual problem which sometimes is easier to solve.
This will be pursued in a later paper.
\end{abstract}

\maketitle

\section{Introduction}\label{sect:intro}
In a real projective space $\p$, a pair of distinct points determines a unique projective line passing through them.
However, there exist two line segments joining them because a real projective line is homeomorphic to a circle.
A natural question arises how we can introduce a notion of convex sets in $\p$.
The first definition probably dates back to Steinitz~\cite{steiniz:conv} and Veblen-Young~\cite{veb:projgeo} according to Danzer, Gr{\"u}nbaum, and Klee~\cite[p.\,159]{danzer:convex} and Deumlich-Elster~\cite{de:conv}.
The following condition is equivalent to that of Steinitz (see Subsection~\ref{subsect:conv}):
a subset $C$ of $\p$ is a {\em (projective) convex set} if, for any two points of $C$, exactly one of the two line segments joining them is contained in $C$.

A convex set $C$ in $\p$ is contained in some affine subspace of $\p$ according to de Groot and de Vries~\cite[Theorem~4]{groot:convex}, which implies that $C$ is a convex set in the usual sense in this affine subspace.
However, for a collection of convex sets $C_{1},\ldots, C_{k}$ in $\p$, there may not exist any affine subspace of $\p$ containing all of them simultaneously.
Also, the intersection of these convex sets may split into multiple connected components each of which is a convex set.
Hence, we call an intersection of a nonempty family of convex sets a {\em multi-convex set}.
Our main interest concerns saturated multi-convex sets, which will be defined in Section~\ref{sect:satmultconv}.

The notion of multi-convex sets appears in the following problem in computational geometry: given a collection of $k$ objects, say polygons, $S_{1},\ldots, S_{k}$ in $\real^{2}$, construct a representation of all lines avoiding every object $S_{i}\ (1\leq i\leq k)$.
By extending $\real^{2}$ to the projective plane $\p$, such a set of lines can be represented as a set of points in the dual projective plane $\pstar$ through the well-known one-to-one correspondence between lines in $\p$ and points in $\pstar$.
Hence, it is useful to introduce the following function $\Phi\colon2^{\p}\to 2^{\pstar}$.
For a subset $S$ of $\p$, we define
\begin{equation*}
\Phi(S) =\set{w^{\ast}\in\pstar}{\text{$w$ is a line in $\p$ avoiding $S$}},
\end{equation*}
where $w^{\ast}$ denotes the point in $\pstar$ corresponding to the line $w$ in $\p$.
Notably $\Phi$ sends a convex set in $\p$ to a convex set in $\pstar$, which we will prove in general (see Proposition~\ref{prop:convtoconv}).
Thus, the problem above is transformed into the problem computing the intersection of $k$ convex sets $\Phi(S_{1}),\ldots,\Phi(S_{k})$ in $\pstar$, which turns out to be a multi-convex set. Moreover, this intersection is saturated as we will see in Section~\ref{sect:satmultconv}.
This problem will be pursued in a later paper.
We remark that our problem is in contrast to problems concerning, for the case of $\real^{2}$, common line transversals, that is, lines meeting every object $S_{i}\ (1\leq i\leq k)$.
Many authors have investigated these problems from computational or combinatorial aspects (see \cite{gruber:convgeo}, \cite[chapter~XIV]{edels:algconv}, \cite{goodmanpollack:convex}, \cite{grun:conv}, \cite{wenger:discrete}, \cite{wenger:transversal}).

In this paper, we study fundamental properties of multi-convex sets and derive a duality for saturated multi-convex sets:
there exists an order anti-isomorphism between nonempty saturated multi-convex sets in $\p$ and those in $\pstar$.
Moreover, we show that the following two antithetical notions are reversed through this duality:
components and co-components of a saturated multi-convex set.
The notions and results studied in this paper will serve as useful tools in studying the problem presented above.

In a complex projective space, the notion of linearly convex sets has been studied in \cite{martin:convex} and \cite{znam:convex}.
We shall show in Remark~\ref{rem:linconv} that the notion of saturated multi-convex sets is almost identical to that of linearly convex sets defined in a real projective space.

In Section~\ref{sect:conv}, we study the order structure of the family of convex sets in $\p$. We focus on the notion of irreducible convex sets, and show that irreducible convex sets are the complements of projective hyperplanes.
In Section~\ref{sect:multconv}, we study multi-convex sets and show two separation properties for multi-convex sets.
In Section~\ref{sect:satmultconv}, we study saturated multi-convex sets and derive the duality for saturated multi-convex sets.

\subsection{Preliminaries}\label{subsect:pre}
In this paper, we only consider finite-dimensional real projective spaces.
Let $V$ be an $(n+1)$-dimensional real vector space.
We denote by $\pn$, $\p$ for short, the $n$-dimensional real projective space associated with $V$, and by $\pi$ the projection of $V\setminus\{0\}$ onto $\pn$, which sends each nonzero vector in $V$ to the linear subspace spanned by the vector.
The dual projective space of $\p$ is denoted by $\pstar$.

\begin{figure}[tb]
  \begin{center}
    \includegraphics[width=4cm]{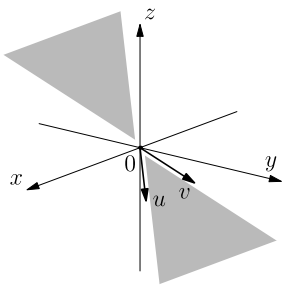}
  \end{center}
  \caption{For two nonzero vectors $u,v\in\real^{3}$, the gray area shows the set of all vectors $\lambda u+\mu v\in\real^{3}$ such that $\lambda,\mu\in\real$, $\lambda\mu\geq 0$, and both of them are not zero.}\label{fig:segment}
\end{figure}

Let $p, q$ be distinct two points in $\p$, and let $u, v$ be nonzero vectors in $V$ such that $\pi(u)=p$ and $\pi(v)=q$.
The {\em (projective) line segments} joining $p$ and $q$ are defined as follows (see Figure~\ref{fig:segment}):
\begin{align*}
&\set{ \pi(\lambda u + \mu v)\in\p}{\text{$\lambda,\mu\in\real$, $\lambda\mu\geq 0$ and ($\lambda\not=0$ or $\mu\not=0$)}};\\
&\set{ \pi(\lambda u + \mu v)\in\p}{\text{$\lambda,\mu\in\real$, $\lambda\mu\leq 0$ and ($\lambda\not=0$ or $\mu\not=0$)}}.
\end{align*}

In an $n$-dimensional affine or projective space, an $(n-1)$-dimensional subspace is called a {\em hyperplane}.
In order to avoid confusing geometric objects in an affine space with those in a projective space, they are sometimes prefixed with the word ``affine'' or ``projective'' such as affine hyperplanes or projective hyperplanes.

Let $(L, \leq)$ and $(M, \leq)$ be partially ordered sets.
A subset $S$ of $L$ is {\em directed} if it is nonempty and every finite subset of $S$ has an upper bound in $S$.
A partially ordered set is a {\em meet-semilattice} if any two elements have a greatest lower bound, that is, an infimum.
A function $f\colon L\rightarrow M$ is called to be {\em order preserving} or {\em monotone} if $x\leq y$ always implies $f(x)\leq f(y)$, and {\em order reversing} or {\em antitone} if $x\leq y$ always implies $f(y)\leq f(x)$.
A function $f\colon L\rightarrow M$ is called a {\em (order) anti-isomorphism} if it is a bijection and both $f$ and $f^{-1}$ are antitone.
A pair $(g,d)$ of functions $g\colon L\rightarrow M$ and $d\colon M\rightarrow L$ is a {\em Galois connection} between $L$ and $M$ if it satisfies the following two conditions:
\begin{enumerate}
\item[(1)] both $g$ and $d$ are monotone;
\item[(2)] the relations $m\leq g(l)$ and $d(m)\leq l$ are equivalent for all pairs of elements $(l,m)\in L\times M$.
\end{enumerate}

\begin{theorem}[{\cite[Theorem\,O-3.6]{gierz:lattice}}]\label{theo:galois}
For every pair of order preserving functions between posets, $g\colon L\to M$ and $d\colon M\to L$, the following conditions are equivalent:
\begin{enumerate}
\item[(i)] $(g,d)$ is a Galois connection;
\item[(ii)] $d\circ g\leq 1_{L}$ and $1_{M} \leq g\circ d$,
\end{enumerate}
where $1_{L}$ and $1_{M}$ denote the identities on $L$ and on $M$.
Moreover, these conditions imply 
\begin{enumerate}
\item[(iii)] $d=d\circ g\circ d$ and $g=g\circ d\circ g$,
\item[(iv)] $g\circ d$ and $d\circ g$ are idempotent.
\end{enumerate}
\end{theorem}

\section{Local algebraic closure systems consisting of projective convex sets}\label{sect:conv}
In this section, we study the order structure of the family of convex sets in $\p$.
We focus on the notion of irreducible convex sets, which is introduced in terms of order, and show that irreducible convex sets are the complements of projective hyperplanes.
This result leads to a one-to-one correspondence between irreducible convex sets in $\p$ and points in $\pstar$.
In Section~\ref{sect:satmultconv}, we will extend this correspondence to a duality for saturated multi-convex sets.

\subsection{Projective convex sets}\label{subsect:conv}
A subset $S$ of $\p$ is {\em semiconvex} if any two points in $S$ can be joined by a line segment which is contained in $S$.
\begin{definition}[Steinitz~\cite{steiniz:conv}]
A subset $S$ of $\p$ is a {\em (projective) convex set} if 
$S$ is semiconvex and there exists a projective hyperplane avoiding $S$.
\end{definition}
Since the complement of a projective hyperplane is an affine subspace of $\p$,
a projective convex set in $\p$ turns out to be a convex set in the usual sense in some affine subspace of $\p$.
The definition of projective convex sets is equivalently described as follows.
\begin{theorem}[de Groot and de Vries~{\cite[Theorem\,4]{groot:convex}}]
A subset $S$ of $\p$ is a projective convex set if and only if, for any two points in $S$, exactly one of the two line segments joining them is contained in $S$.
\end{theorem}
According to this theorem, the notion of projective convex sets can be characterized by only the notion of line segments: we need not refer to hyperplanes in order to define convex sets.
We could adopt it as the definition of projective convex sets.

\subsection{Local algebraic closure systems}
Let $X$ be a set, and let $\mathcal{S}$ be a family of subsets of $X$.
The family $\mathcal{S}$ is called an {\em algebraic closure system} on $X$ if it has the following three properties:
\begin{enumerate}
\item[(C1)] the set $X$ is in $\mathcal{S}$;
\item[(C2)] the intersection of any nonempty subfamily of $\mathcal{S}$ is in $\mathcal{S}$;
\item[(C3)] the union of any directed subfamily of $\mathcal{S}$ is in $\mathcal{S}$.
\end{enumerate}
An algebraic closure system is referred to as an {\em alignment}, a {\em topped algebraic $\bigcap$-structure} , or a {\em convexity} by some authors (see Coppel~\cite[chapter~I]{coppel:convgeo}, Davey-Priestley~\cite[p.\,150]{davey:lattice-order}, and Van de Vel~\cite{vandevel:convstruct}).
The family of convex sets in $\rn$ is a typical example of an algebraic closure system.
An algebraic closure system is known to provide various notions such as extreme point, independent set, basis and face.

We have to take into account that projective convex sets do not form an algebraic closure system.
Hence, it is worth considering what kind of set-system they form.

\begin{example}\label{exam:notalgsys}
The family of convex sets in $\p$ is not an algebraic closure system on $\p$.
\begin{enumerate}
\item Clearly the whole space $\p$ is not convex.
\item For two distinct points in $\p$, the intersection of the two line segments joining them is the two-point set. Thus, the intersection of two convex sets need not be convex.
\item In $\p^{1}=\real\cup\{\infty\}$, the union of convex sets $\p^{1}\setminus \left(0,\tfrac{1}{n}\right)$ for $n\in\nat$ covers the whole space $\p^{1}$, where $\left(0,\tfrac{1}{n}\right)$ is an open interval in $\real$.
Thus, the union of a directed family of convex sets need not be convex.
\end{enumerate}
\end{example}

The following proposition states that the two conditions (C2) and (C3) above can be satisfied if a collection $\mathcal{F}$ of convex sets in $\p$ is consistent, where $\mathcal{F}$ is called to be {\em consistent} if there exists an affine subspace of $\p$ containing all members of $\mathcal{F}$ simultaneously.

\begin{proposition}\label{prop:localgsys}
Let $\mathcal{S}$ be the family of all convex sets in $\p$. Then we have the following two properties:
\begin{enumerate}
\item the intersection of any nonempty consistent subfamily of $\mathcal{S}$ is in $\mathcal{S}$;
\item the union of any consistent directed subfamily of $\mathcal{S}$ is in $\mathcal{S}$.
\end{enumerate}
\end{proposition}

One can say that the family of convex sets in $\p$ is a local algebraic closure system on $\p$:
the whole space $\p$ is covered by open affine subspaces, and ``locally'', that is, on each of these open affine subspaces, the convex sets form an algebraic closure system.
This observation leads to the following notion of local algebraic closure systems:

\begin{definition}
Let $X$ be a set, and let $\mathcal{S}$ be a family of subsets of $X$.
Then $\mathcal{S}$ is called a {\em local algebraic closure system} on $X$ if it satisfies the two conditions in Proposition~\ref{prop:localgsys}, where a collection $\mathcal{F}$ of $\mathcal{S}$ is called to be {\em consistent} if there exists a member of $\mathcal{S}$ containing all members of $\mathcal{F}$ simultaneously.
\end{definition}
Clearly an algebraic closure system is a local algebraic closure system, however the converse is not true in general as Example~\ref{exam:notalgsys} shows.

We remark that the hypothesis of consistency in Proposition~\ref{prop:localgsys}(ii) can be dropped for open convex sets:

\begin{proposition}\label{prop:directclosed}
The union of any directed family of open convex sets is an open convex set.
\end{proposition}

\begin{proof}
Let $\mathcal{F}$ be a directed family of open convex sets.
For two points $p, q\in\bigcup\mathcal{F}$, let $U$ and $U'$ be members of $\mathcal{F}$ containing $p$ and $q$, respectively.
Since $\mathcal{F}$ is directed, there exists a member $V$ of $\mathcal{F}$ containing both $U$ and $U'$.
Then $V$ contains one of the two line segments joining $p$ and $q$.
This line segment is contained in $\bigcup\mathcal{F}$.

To prove the uniqueness, assume that $\bigcup\mathcal{F}$ contains both of the two line segments joining $p$ and $q$.
Then $\bigcup\mathcal{F}$ contains the whole line spanned by $p$ and $q$.
By compactness there exists a member of $\mathcal{F}$ containing this whole line.
However, this contradicts that every member of $\mathcal{F}$ is convex.
We have proved the proposition.
\end{proof}

\subsection{Irreducible convex sets}\label{subsect:irred}
The notion of irreducible convex sets we introduce now appears in studying multi-convex sets:
when we have a multi-convex set $C$, we want to associate it with a family $\mathcal{S}$ of convex sets satisfying $C=\bigcap\mathcal{S}$;
this family $\mathcal{S}$ of $C$ is reduced to another one if some element $N$ in $\mathcal{S}$ can separate into a pair of strictly greater convex sets $N_{1}$ and $N_{2}$ satisfying $N=N_{1}\cap N_{2}$;
if this reduction terminates, then we have such convex sets that can not be the intersection of any pair of strictly greater convex sets.
One can say that convex sets of this kind are irreducible.

\begin{definition}
A convex set $A$ is {\em irreducible} if, for any convex sets $B$ and $C$, $A=B\cap C$ implies $A=B$ or $A=C$.
\end{definition}

In a meet-semilattice, an element $a$ is called to be {\em irreducible} if $a=b\wedge c$ always implies $a=b$ or $a=c$ (see \cite[p.\,53]{davey:lattice-order} and \cite[section~I-3]{gierz:lattice}).

\begin{theorem}\label{theo:irred}
For a subset $C$ of $\p$, the following conditions are equivalent:
\begin{enumerate}
\item $C$ is an irreducible convex set;
\item $C$ is a maximal convex set;
\item $C$ is the complement of a projective hyperplane.
\end{enumerate}
\end{theorem}

\begin{proof}
That (iii) $\Rightarrow$ (ii) $\Rightarrow$ (i) is trivial.
We prove that (i) implies (iii).
Assume that $C$ is an irreducible convex set, and not the complement of any projective hyperplane.
Let $A$ be a maximal affine subspace, which is the complement of a hyperplane, of $\p$ which contains $C$.
Since $C$ is a proper subset of $A$, clearly the closure of $C$ relative to $A$ is not $A$.
Let $p$ be a point in $A$ which is not adherent to $C$.
It is known that a closed affine convex set is the intersection of the closed half-spaces which contain it (see \cite[Theorem 11.5]{rock:convanly}).
Thus we can easily derive that there exists an affine hyperplane in $A$ passing through $p$ and avoiding $C$.
By extending such an affine hyperplane to the projective hyperplane, we have another affine subspace $A'$ of $\p$ containing $C$, which avoids $p$ by construction.

Let us choose an arbitrary point $o$ in $C$, and let $u$ be an arbitrary nonzero vector in $A'$ emanating from $o$ to $p$.
We define the following two cones in $A'$:
\begin{align*}
{\uparrow}C   & =\set{x+\lambda u\in A'}{\text{$x\in C$ and $\lambda\geq 0$}};\\
{\downarrow}C & =\set{x-\lambda u\in A'}{\text{$x\in C$ and $\lambda\geq 0$}}.
\end{align*}
Trivially we have $C\subseteq{\uparrow}C \cap{\downarrow}C$.
We prove $C\supseteq {\uparrow}C \cap{\downarrow}C$.
Any point $z$ in ${\uparrow}C \cap{\downarrow}C$ has the following two forms: $x+\lambda u$ and $y-\lambda' u$ for $x,y\in C$, $\lambda\geq 0,$ and $\lambda'\geq 0$.
Hence, we have $y=x+(\lambda+\lambda')u$.
This implies that the point $z$ lies in the line segment joining $x$ and $y$ which is contained in $C$.
Thus we have $C={\uparrow}C \cap{\downarrow}C$.
Since $p$ is not adherent to $C$, both ${\uparrow}C$ and ${\downarrow}C$ are strictly greater than $C$.
However, these results contradict that $C$ is irreducible.
Therefore, if $C$ is irreducible, then it is the complement of a projective hyperplane.
\end{proof}

\begin{remark}\label{rem:irred}
The irreducibility of convex sets can be reduced to the irreducibility of open convex sets, where 
an open convex set $A$ is called to be {\em irreducible} if, for any open convex set $B$ and $C$, $A=B\cap C$ implies $A=B$ or $A=C$.
Clearly any irreducible convex set is an irreducible open convex set by definition.
The converse can be easily proved in a similar way to the proof of Theorem~\ref{theo:irred}.
Thus it follows that a subset of $\p$ is an irreducible open convex set if and only if it is an irreducible convex set.
Therefore, in order to prove that a given convex set $S$ is irreducible, 
it suffices to prove that $S$ is an irreducible open convex set.
\end{remark}

The projective duality between points and hyperplanes states that there exists a natural one-to-one correspondence between points in $\p$ and hyperplanes in $\pstar$ (and equivalently between hyperplanes in $\p$ and points in $\pstar$) which reverses the incidence relation: a point $p$ in $\p$ lies in a hyperplane $h$ in $\p$ if and only if the point in $\pstar$ corresponding to $h$ lies in the hyperplane in $\pstar$ corresponding to $p$.
From this correspondence and Theorem~\ref{theo:irred}, we immediately obtain a one-to-one correspondence between points in $\p$ and irreducible convex sets in $\pstar$ (and between irreducible convex sets in $\p$ and points in $\pstar$).

\begin{notation}
When there is no danger of confusion, we denote by $\delta(p)$ and $\delta(C)$ the irreducible convex set and the point in $\pstar$ corresponding to a point $p$ and an irreducible convex set $C$ in $\p$, respectively.
Similarly we denote by $\delta(p')$ and $\delta(C')$ the irreducible convex set and the point in $\p$ corresponding to a point $p'$ and an irreducible convex set $C'$ in $\pstar$, respectively.
\end{notation}

\begin{proposition}\label{prop:projduality}
For a point $p$ in $\p$ and for an irreducible convex set $C$ in $\p$, we have the following two properties:
\begin{enumerate}
\item $\delta(\delta(p))=p$ and $\delta(\delta(C))=C$;
\item the relations $p\in C$ and $\delta(C)\in \delta(p)$ are equivalent.
\end{enumerate}
\end{proposition}

This result will be extended to a duality for saturated multi-convex sets in Section~\ref{sect:satmultconv}.

\subsection{Convex hulls relative to convex sets}
A subset $S$ of an affine space always has a least convex set containing $S$, which is called the convex hull of $S$, while a subset $S$ of $\p$ need not: for example, distinct two points are connected by distinct two line segments.
However, ``locally'', that is, on each convex set containing $S$, we have the following notion of convex hull of $S$.

\begin{definition}
Let $C$ be a convex set in $\p$.
For a subset $S$ of $C$, the {\em convex hull} $\relhull{S}{C}$ of $S$ {\em relative to} $C$ is a least convex subset of $C$ containing $S$:
\begin{equation*}
\relhull{S}{C}=\bigcap\set{N}{\text{$N$ is a convex set satisfying $S\subseteq N\subseteq C$}}.
\end{equation*}
\end{definition}

In particular, for two points $p$ and $q$ in a convex set $C$, we denote by $\relhull{p,q}{C}$ the line segment joining $p$ and $q$ which is contained in $C$.

\begin{proposition}
Let $C$ be a convex set in $\p$. For subsets $S$ and $T$ of $C$, we have the following four properties:
\begin{enumerate}
\item $S\subseteq \relhull{S}{C}$;
\item $S\subseteq T$ implies $\relhull{S}{C}\subseteq \relhull{T}{C}$;
\item $\relhull{\relhull{S}{C}}{C}=\relhull{S}{C}$;
\item $\relhull{S}{C}=\bigcup\set{\relhull{F}{C}}{\text{$F\subseteq S$ and $F$ is finite}}$.
\end{enumerate}
\end{proposition}

\begin{proof}
By Proposition~\ref{prop:localgsys}, the family $\mathcal{S}$ of all convex subsets of $C$ is an algebraic closure system on $C$.
The first two properties (i) and (ii) are trivial.

(iii) It follows from (C1) and (C2) that $\relhull{S}{C}$ is in $\mathcal{S}$.
Hence we obtain $\relhull{\relhull{S}{C}}{C}=\relhull{S}{C}$.

(iv) Let $\mathcal{T}$ be the family of all convex sets $\relhull{F}{C}$ such that $F$ is a finite subset of $S$.
Clearly we have $\relhull{S}{C}\supseteq\bigcup\mathcal{T}$.
Since $S\subseteq\bigcup\mathcal{T}$, we derive $\relhull{S}{C}\subseteq \relhull{\bigcup\mathcal{T}}{C}$ from the property (ii).
Since $\mathcal{T}$ is directed, its union is in $\mathcal{S}$ by (C3), which implies $\relhull{\bigcup\mathcal{T}}{C}=\bigcup\mathcal{T}$.
Thus we have $\relhull{S}{C}=\bigcup\mathcal{T}$.
\end{proof}

We remark that the proposition above is derived from only the fact that convex sets in $\p$ form a local algebraic closure system.

In an affine space, the convex hull of an open set is open, and that of a compact set is compact, both of which are essentially due to Carath{\'e}odory's theorem (see \cite[section~17]{rock:convanly}), although the convex hull of a closed set need not be closed: for example, the convex hull of $S$ is not closed when $S$ is the union of a closed half-space and a single point not on the half-space.
Since $\p$ is a compact Hausdorff space, a closed set in $\p$ is compact, and vice versa.
Thus we have:

\begin{proposition}\label{prop:convhull}
Let $C$ be a convex set in $\p$. For a subset $S$ of $C$, we have the following two properties:
\begin{enumerate} 
\item if $S$ is closed in $\p$, then the convex hull of $S$ relative to $C$ is closed in $\p$;
\item if $S$ is open in $\p$, then the convex hull of $S$ relative to $C$ is open in $\p$.
\end{enumerate}
\end{proposition}

\section{Multi-convex sets and their separation properties}\label{sect:multconv}
In this section, we introduce two antithetical notions: components and co-components of a multi-convex set.
We show two separation properties for multi-convex sets, which motivate the notion of saturated multi-convex sets in Section~\ref{sect:satmultconv}.

\subsection{Multi-convex sets}
\begin{definition}
A subset of $\p$ is a {\em multi-convex set} if it is the intersection of a nonempty family of projective convex sets.
\end{definition}

\begin{example}
Figure~\ref{fig:convcomp-c} shows an example of a disconnected multi-convex set in $\p^{2}$.
Let $C$ and $D$ be two convex sets in $\p^{2}$ such that there exists no convex set containing both of them as in Figure~\ref{fig:convcomp-a} and~\ref{fig:convcomp-b}.
By definition, there exist two hyperplanes which avoid $C$ and $D$, respectively.
The intersection of $C$ and $D$ consists of two connected regions which are separated by two hyperplanes.
\begin{figure}[tb]
  \begin{minipage}{0.32\hsize}
    \begin{center}
      \subfigure[]{\includegraphics[height=3.5cm]{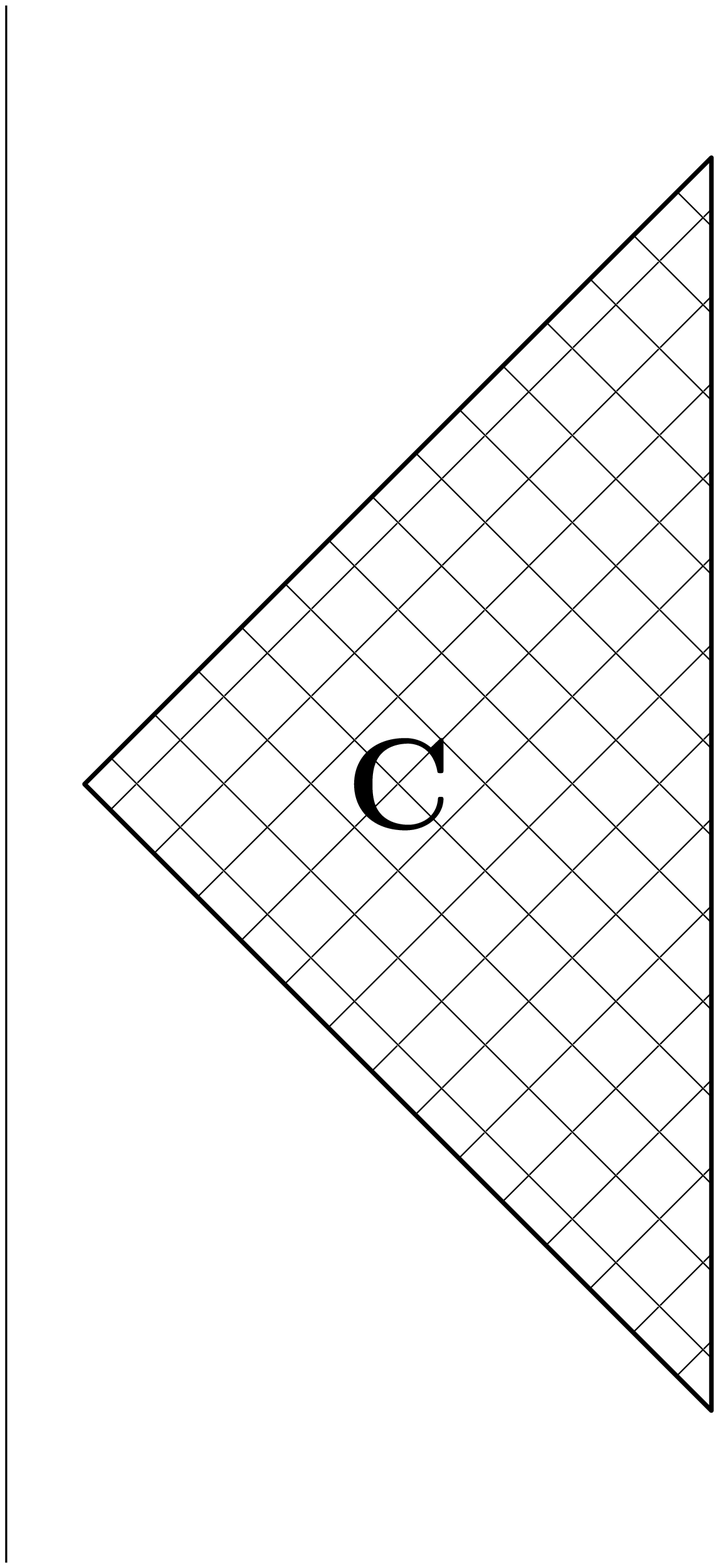}\label{fig:convcomp-a}}
    \end{center}
  \end{minipage}
  \begin{minipage}{0.32\hsize}
    \begin{center}
      \subfigure[]{\includegraphics[height=3.5cm]{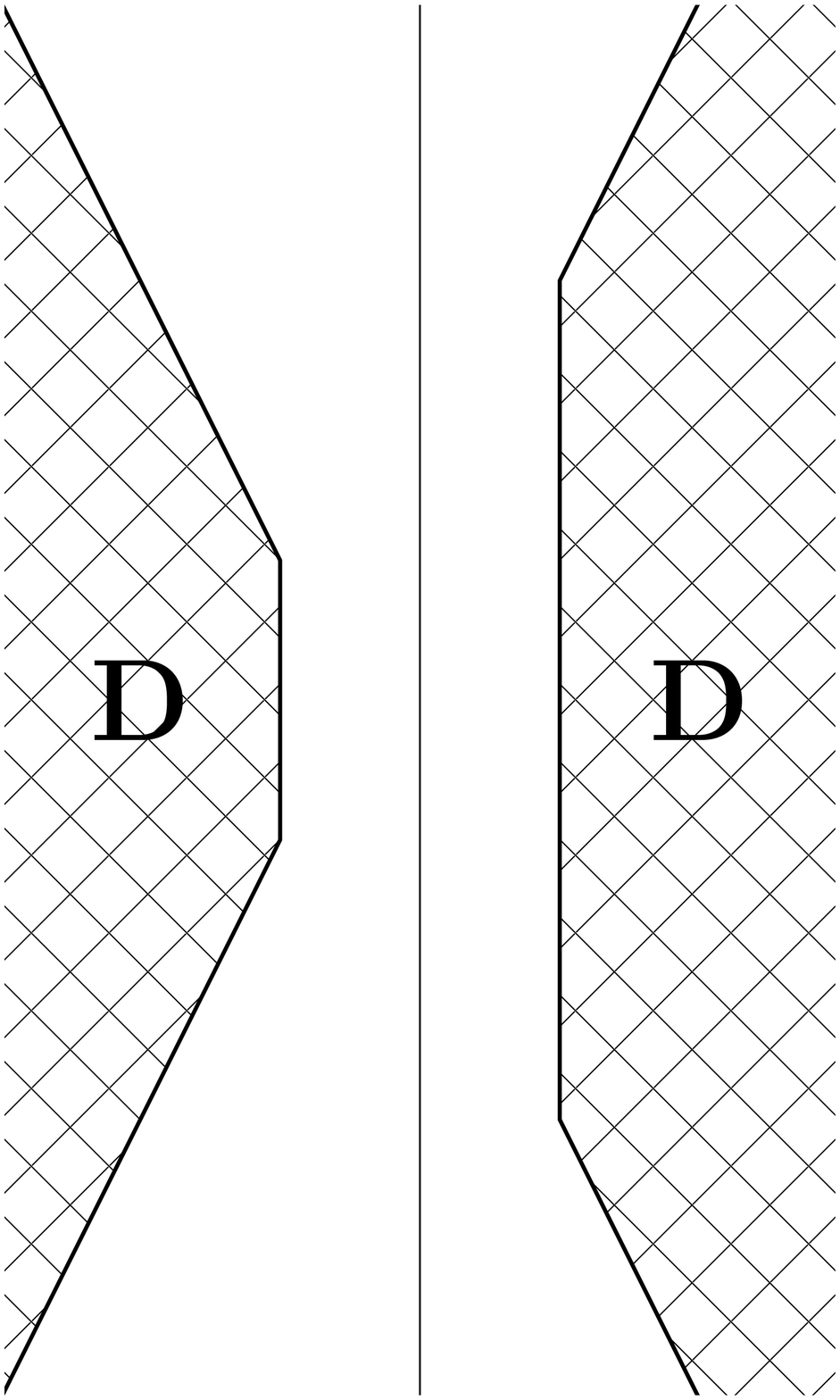}\label{fig:convcomp-b}}
    \end{center}
  \end{minipage}
  \begin{minipage}{0.32\hsize}
    \begin{center}
      \subfigure[]{\includegraphics[height=3.5cm]{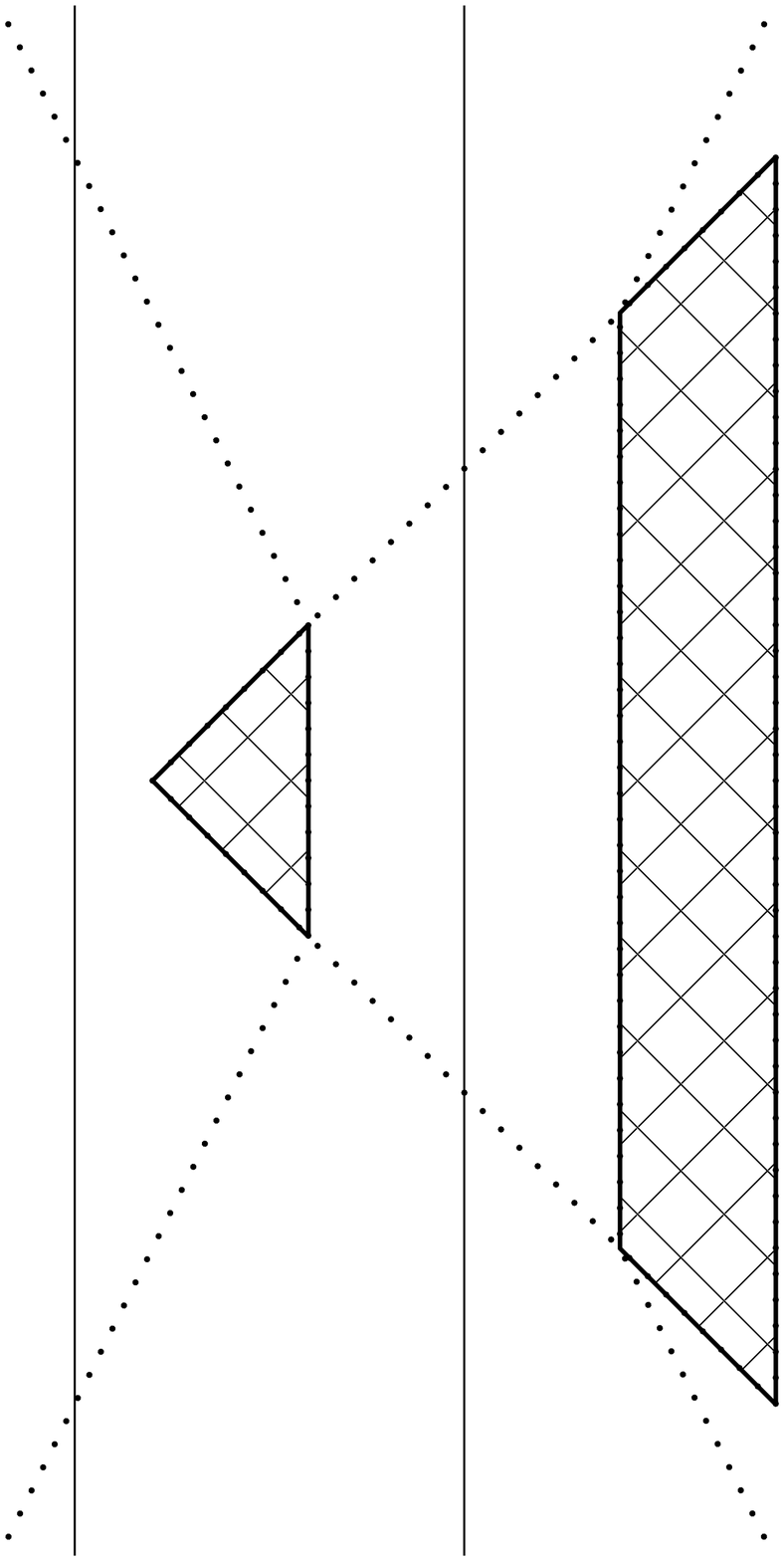}\label{fig:convcomp-c}}
    \end{center}
  \end{minipage}
  \caption{(a) A convex set $C$ with a line avoiding $C$, (b) a convex set $D$ with a line avoiding $D$, and (c) the intersection of $C$ and $D$ whose two connected components are separated by these two lines.}\label{fig:convcomp}
\end{figure}
\end{example}

\begin{proposition}\label{prop:comp}
Every connected component of a multi-convex set is a convex set.
\end{proposition}

\begin{proof}
Let $C$ be a multi-convex set, and let $\{C_{i}\}_{i\in I}$ be a family of convex sets such that $C=\bigcap_{i\in I}C_{i}$.
For any two points $p,q\in C$, if there exist two members $C_{i}$ and $C_{j}\ (i,j\in I)$ such that $\relhull{p,q}{C_{i}}\not=\relhull{p,q}{C_{j}}$,
then $p$ and $q$ lie in different connected components of $C_{i}\cap C_{j}$, hence $p$ and $q$ line in different connected components of $C$.
It follows that if $p$ and $q$ lie in the same connected component $N$ of $C$, all the segments $\relhull{p,q}{C_{i}}\ (i\in I)$ coincide. This implies that the line segment is contained in $N$, joining $p$ and $q$.
Clearly $N$ does not contain both of the two line segments joining $p$ and $q$.
We have proved the proposition.
\end{proof}

Note that every connected component of a multi-convex set $C$ is a maximal convex subset of $C$, and vice versa.
We introduce the following notion of components in terms of order rather than topology: this view allows to dualize the notion of components as we shall see later.

\begin{definition}
Let $C$ be a multi-convex set.
We call every maximal convex subset of $C$ a {\em component} of $C$ and denote by $\comp{C}$ the family of all components of $C$.
We call the cardinality of $\comp{C}$ the {\em degree} of $C$.
\end{definition}

From Proposition~\ref{prop:comp}, we immediately derive the following result.

\begin{corollary}
Let $C$ be a multi-convex set.
Then $\comp{C}$ is a mutually disjoint family of convex sets such that $C=\bigcup\comp{C}$.
\end{corollary}

\begin{remark}
We have to take into account that a disjoint union of convex sets need not be a multi-convex set (see Figure~\ref{fig:conv-4deg}).
Let us consider three lines in general position in $\p^{2}$.
The complement of the union of these three lines is a multi-convex set of degree $4$, however the disjoint union $S$ of all except for the $4$th component is not a multi-convex set.
Because any convex set containing $S$ is the complement of one of those three lines, and the least multi-convex set containing $S$ must contain the $4$th component.
\begin{figure}[tb]
  \begin{minipage}{0.32\hsize}
    \begin{center}
      \subfigure[]{\includegraphics[height=3.5cm]{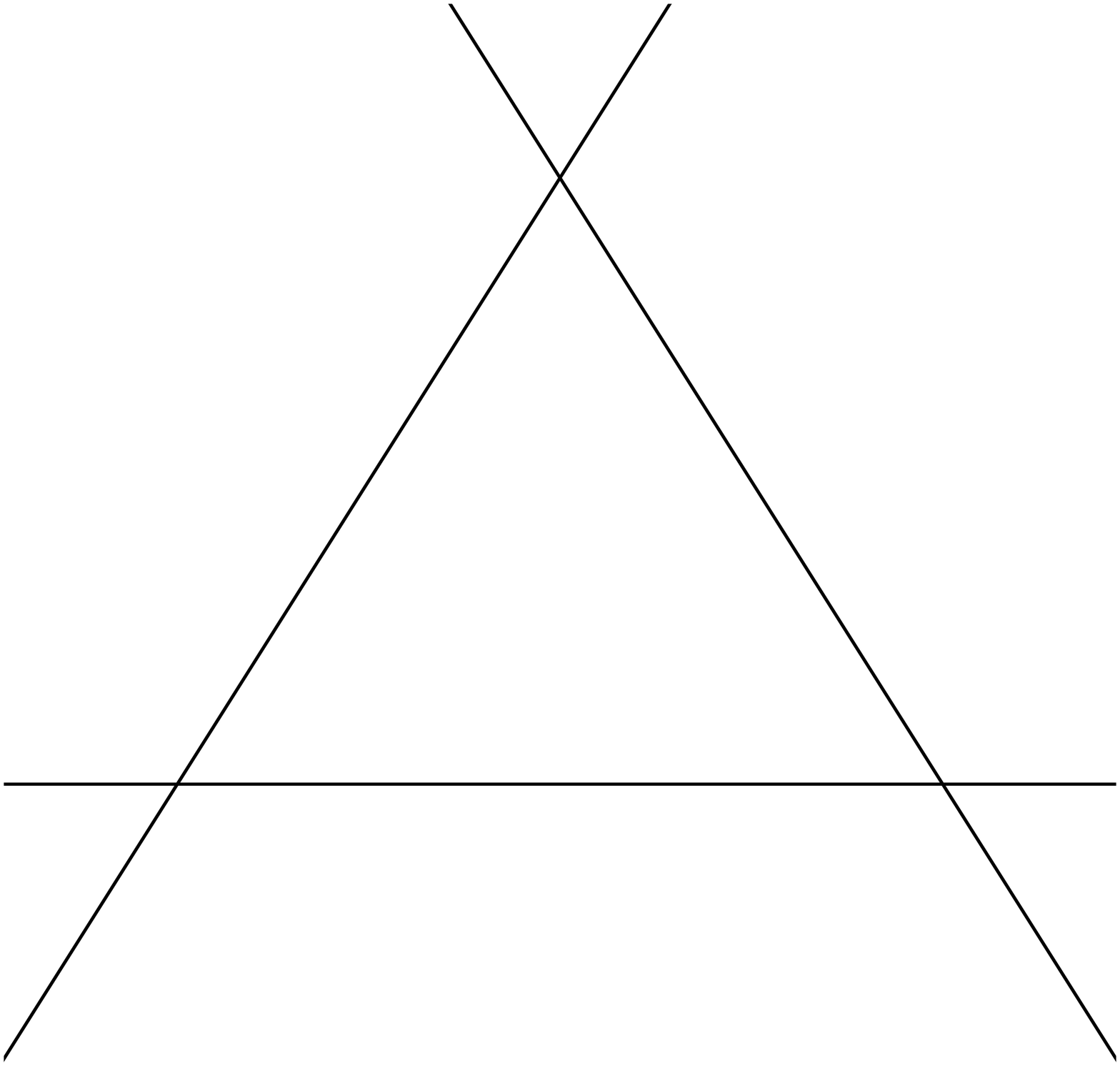}\label{fig:conv-4deg-a}}
    \end{center}
  \end{minipage}
  \begin{minipage}{0.32\hsize}
    \begin{center}
      \subfigure[]{\includegraphics[height=3.5cm]{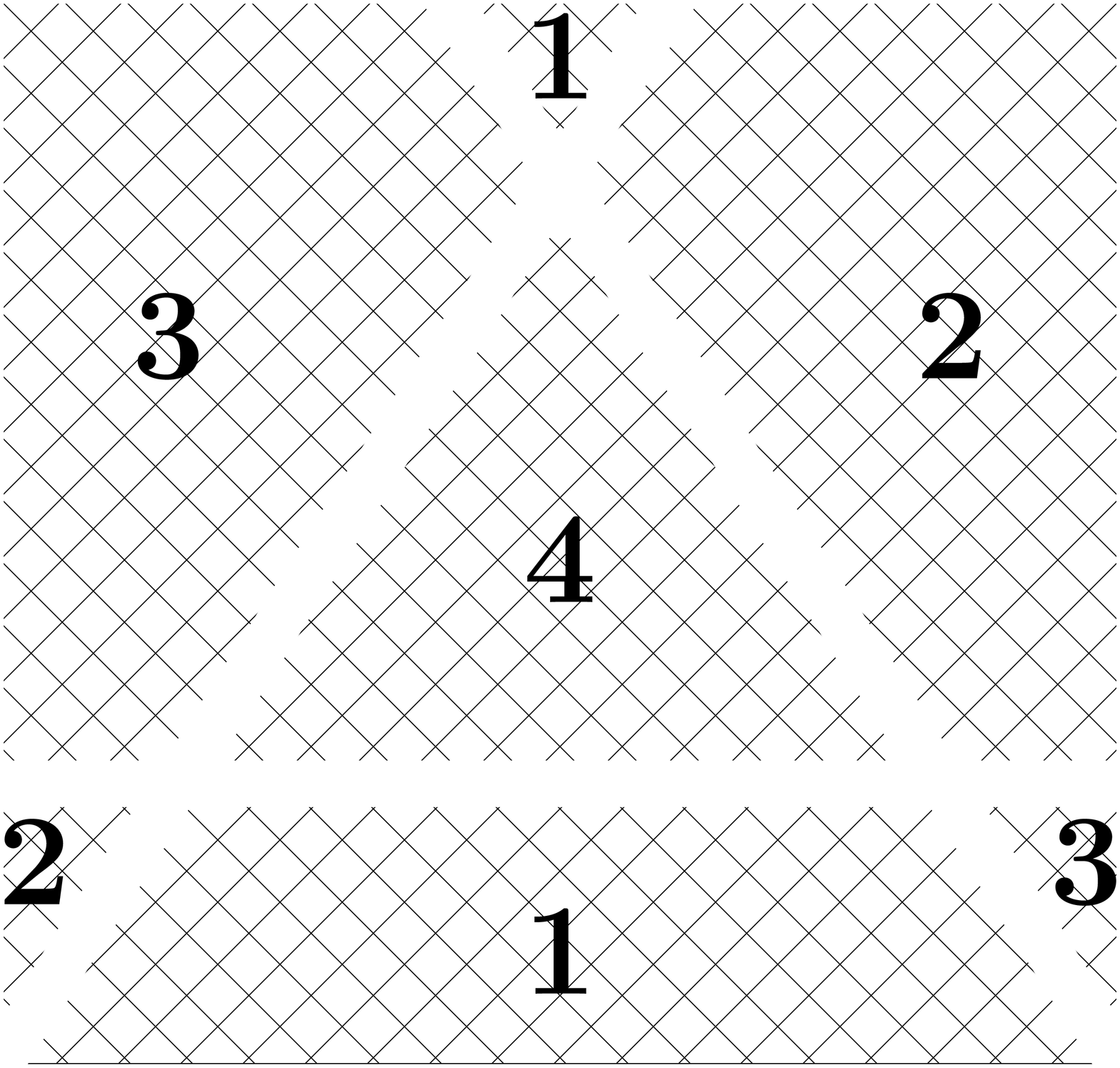}\label{fig:conv-4deg-b}}
    \end{center}
  \end{minipage}
  \begin{minipage}{0.32\hsize}
    \begin{center}
      \subfigure[]{\includegraphics[height=3.5cm]{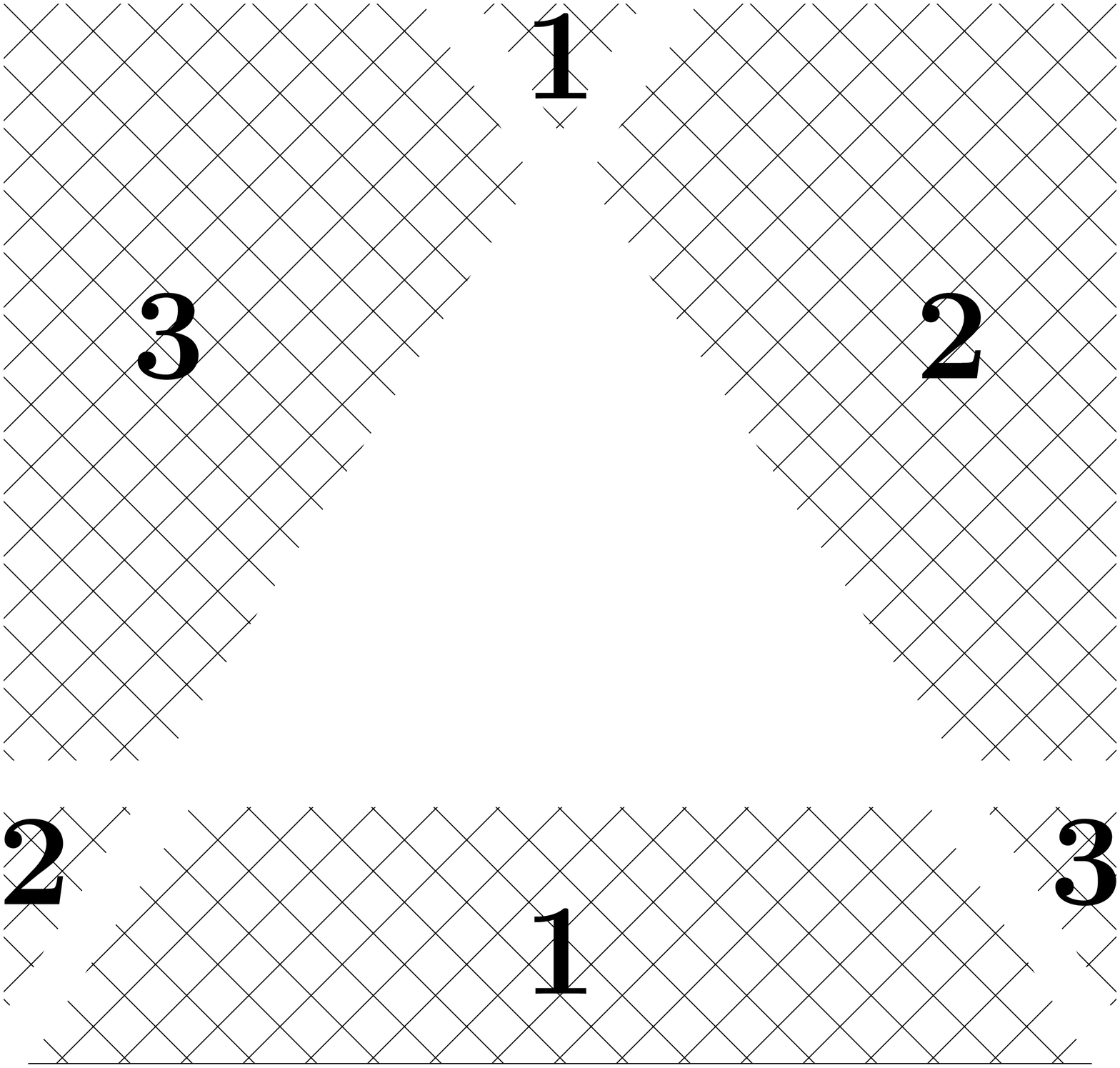}\label{fig:conv-4deg-c}}
    \end{center}
  \end{minipage}
  \caption{(a) Three projective lines in general position, (b) the complement of the union of these three lines, and (c) the disjoint union of all except for the $4$th component.}\label{fig:conv-4deg}
\end{figure}
\end{remark}

\begin{definition}
Let $C$ be a multi-convex set.
We call every minimal convex set containing $C$ a {\em co-component} of $C$ and denote by $\cocomp{C}$ the family of all co-components of $C$.
We call the cardinality of $\cocomp{C}$ the {\em co-degree} of $C$.
\end{definition}

Figure~\ref{fig:convcomp} shows that $C$ and $D$ are the co-components of the multi-convex set $C\cap D$.

A collection of convex sets is called to be {\em mutually inconsistent} if there is no convex set in $\p$ containing any pair of convex sets in this collection.

\begin{proposition}\label{prop:cocomp}
Let $C$ be a multi-convex set.
Then $\cocomp{C}$ is a mutually inconsistent family of convex sets such that $C=\bigcap\cocomp{C}$.
\end{proposition}

\begin{proof}
For two members $N,N'\in\cocomp{C}$, if $N$ and $N'$ are contained in some convex set $M$, then the convex hull $\relhull{C}{M}$ of $C$ relative to $C$ coincides with both $N$ and $N'$, hence we obtain $N=N'$.
Therefore every pair of distinct members in $\cocomp{C}$ is inconsistent.
Let $\mathcal{S}$ be a family of convex sets such that $C=\bigcap\mathcal{S}$.
Clearly we have $\bigcap\cocomp{C}\subseteq\bigcap\mathcal{S}$. Conversely we have $C\subseteq\bigcap\cocomp{C}$ by the definition of co-components.
We have proved $C=\bigcap\cocomp{C}$.
\end{proof}

\subsection{Separation properties for multi-convex sets}

\begin{proposition}\label{prop:separation1}
Let $C$ be a multi-convex set, and let $p$ be a point not in the topological interior $\introp{C}$ of $C$.
Then there exists an open convex set which contains $\introp{C}$ and avoids $p$.
\end{proposition}

\begin{proof}
The topology of $\p$ has a countable base of open convex sets.
Hence, there exists a countable family $\{U_{i}\}_{i\in\nat}$ of open convex sets such that $ \introp{C}=\bigcup_{i\in\nat} U_{i}$.
We introduce the following notation: $C_{k}$ denotes $\bigcup_{i\leq k}U_{i}$, and $\hull{C_{k}}$ denotes the least multi-convex set containing $C_{k}$, that is, the intersection of all convex sets containing $C_{k}$.

We first prove that there exists a co-component of $\hull{C_{k}}$ which avoids $p$.
Since $C_{k}$ is open, every minimal convex set containing $C_{k}$ is open by Proposition~\ref{prop:convhull}.
Since $C_{k}$ consists of finitely many connected components,
the number of all minimal convex sets containing $C_{k}$ is finite.
It follows that $\hull{C_{k}}$ is an open multi-convex set of finite co-degree.
Since $\hull{C_{k}}$ is an open subset of $C$, we obtain $\hull{C_{k}}\subseteq\introp{C}$.
Since $\hull{C_{k}}$ is the intersection of its co-components, there exists a co-components of $\hull{C_{k}}$ which avoids $p$.

We show that, for each $k\in\nat$, we can choose a co-component $N_{k}$ of $\hull{C_{k}}$ avoiding $p$ in such a way that $N_{1}\subseteq N_{2}\subseteq \cdots$.
Let us consider the following directed graph $T$:
the vertices are the pairs $\pair{k}{N}$ , where $k\in\nat$ and $N$ is a co-component of $\hull{C_{k}}$ avoiding $p$; there exists an edge from $\pair{k}{N}$ to $\pair{k'}{N'}$ if $k'=k+1$ and $N\subseteq N'$.
Then $T$ is a finitely branching, infinite, rooted tree.
In deed, as we have proved above, for each $k\in \nat$, there exists a co-component $N$ of $\hull{C_{k}}$ avoiding $p$, hence $\pair{k}{N}$ is a vertex of $T$.
For each $i\leq k$, the convex hull $\relhull{C_{i}}{N}$ of $C_{i}$ relative to $N$ is a co-component of $\hull{C_{i}}$ avoiding $p$, hence $\pair{i}{\relhull{C_{i}}{N}}$ is also a vertex of $T$.
Thus we obtain the following decreasing sequence of vertices starting at $\pair{k}{N}$:
\[
\pair{k}{N}\leftarrow\pair{k-1}{\relhull{C_{k-1}}{N}}\leftarrow\cdots\leftarrow \pair{1}{\relhull{C_{1}}{N}}.
\]
Since the co-component $\relhull{C_{1}}{N}$ of $C_{1}$ coincides with $U_{1}$ by definition, this sequence ends with the root of $T$.
It follows that $T$ is a connected, infinite, rooted graph.
Since $\hull{C_{k}}$ is of finite co-degree for all $k\in\nat$, the graph $T$ is finitely branching.
To prove that $T$ is a tree, assume that there exists a cycle $\mathcal{C}$ in $T$.
Let $\pair{k}{N}$ be a vertex in $\mathcal{C}$ whose index $k$ is maximum among those of all vertices in $\mathcal{C}$.
Then we have distinct two vertices $\pair{k-1}{M}$ and $\pair{k-1}{M'}$ in $\mathcal{C}$ adjacent to $\pair{k}{N}$ such that $\pair{k-1}{M}\rightarrow\pair{k}{N}\leftarrow\pair{k-1}{M'}$.
Since $N$ contains both $M$ and $M'$, we obtain $M=M'$, which is a contradiction.
We have proved that $T$ is a finitely branching, infinite, rooted tree.
K\"{o}nig's lemma states that if a finitely branching rooted tree is infinite, then there exists an infinitely long rooted path (see~\cite[section~3]{hansen:konig}).
Thus we have one.
Clearly such a path is increasing:
\[
\pair{1}{N_{1}}\rightarrow\cdots\rightarrow\pair{k}{N_{k}}\rightarrow\pair{k+1}{N_{k+1}}\rightarrow\cdots.
\]
The family of all the co-components $N_{i}$ is totally ordered and consists of open convex sets avoiding $p$.
Hence the union $\bigcup_{i\in\nat} N_{i}$ is an open convex set avoiding $p$.
Moreover, since $C_{i}\subseteq N_{i}$ for all $i\in\nat$, we obtain $\introp{C}=\bigcup_{i\in\nat} C_{i}\subseteq\bigcup_{i\in\nat}N_{i}$.
\end{proof}

\begin{proposition}\label{prop:separation2}
Let $K$ and $L$ be closed multi-convex sets.
Then the following conditions are equivalent:
\begin{enumerate}
\item there exist a co-component of $K$ and a co-component of $L$ which are disjoint;
\item there exist disjoint open convex sets $U$ and $V$ such that $K\subseteq U$ and $L\subseteq V$.
\end{enumerate}
\end{proposition}

\begin{proof}
Assume the condition (i).
Let $(N,M)$ be a disjoint pair of a co-component of $K$ and that of $L$.
Then there exist disjoint two open sets $S$ and $T$ such that $N\subseteq S$ and $M\subseteq T$.
Let $A$ and $B$ be two affine subspaces of $\p$ containing $N$ and $M$, respectively.
Without loss of generality, we can assume $S\subseteq A$ and $T\subseteq B$.
Let us denote by $V_{A}$ the vector space associating to $A$ .
From the compactness of $N$, we can derive that there exists an open convex subset $W$ of $V_{A}$ such that $N+W\subseteq S$, where $N+W:=\set{p+u}{\text{$p\in N$ and $u\in W$}}$.
Then the open convex set $U:=N+W$ satisfies $N\subseteq U\subseteq S$.
Similarly, we have an open convex set $V$ in $\p$ such that $M\subseteq V\subseteq T$.
Clearly $U$ and $V$ are disjoint.

Conversely, assume the condition (ii).
Then $\relhull{K}{U}$ and $\relhull{L}{V}$ are co-components of $K$ and $L$, respectively.
Since $U$ and $V$ are disjoint, $\relhull{K}{U}$ and $\relhull{L}{V}$ are disjoint.
\end{proof}

\section{Saturated multi-convex sets and their duality}\label{sect:satmultconv}
In this section, we introduce the notion of saturated multi-convex sets, and show that a subset $S$ of $\p$ is a saturated multi-convex set if and only if $S$ is the intersection of a nonempty family of irreducible convex sets containing $S$.
We derive a duality for saturated multi-convex sets: there exists an order anti-isomorphism between nonempty saturated multi-convex sets in $\p$ and those in $\pstar$. Moreover, we show that the notion of components and that of co-components are reversed through this duality.

\subsection{Saturated multi-convex sets}
We introduce a class of multi-convex sets with the following separation property:

\begin{definition}
A multi-convex set $S$ of $\p$ is {\em saturated} if, for each point $p$ not in $S$, there exists an open convex set containing $S$ and avoiding $p$.
\end{definition}
Note that a subset of $\p$ is a saturated multi-convex set if and only if it is the intersection of a nonempty family of open convex sets.

From the two separation properties in the previous section, we can immediately deduce the following two corollaries.

\begin{corollary}
The interior of a multi-convex set is a saturated multi-convex set.
\end{corollary}

Note that the closure of a multi-convex set need not even be a multi-convex set: for example, the closure of the intersection of finitely many irreducible convex sets is the whole space.

\begin{corollary}
An open or closed multi-convex set is saturated.
\end{corollary}

In Subsection~\ref{subsect:comp-cocomp}, we will see examples of multi-convex sets which are not saturated (see Figure~\ref{fig:antiexam-c} and~\ref{fig:antiexam-d}).

In Subsection~\ref{subsect:irred}, we mentioned the reduction of families $\mathcal{S}$ of convex sets satisfying $C=\bigcap\mathcal{S}$ for a multi-convex set $C$.
We now want to know when a multi-convex set is the intersection of such a family of convex sets that can not be reduced any more, that is, the intersection of irreducible convex sets.
In deed, it is saturated one as we show now:

\begin{lemma}\label{lem:separation}
Let $C$ be a saturated multi-convex set, and let $p $ be a point not in $C$.
Then there exists an irreducible convex set which contains $C$ and avoids $p$.
\end{lemma}

\begin{proof}
Let $\mathcal{F}$ be the family of all open convex sets that contains $C$ and avoids $p$.
By Proposition~\ref{prop:directclosed}, the union of any nonempty totally ordered subfamily of $\mathcal{F}$ is an open convex set.
Clearly this open convex set avoids $p$, hence it belongs to $\mathcal{F}$.
By Zorn's lemma, there exists a maximal member $A$ in $\mathcal{F}$.
We prove that $A$ is an irreducible convex set.
By Remark~\ref{rem:irred}, it suffices to prove that $A$ is an irreducible open convex set.
Assume the opposite.
Then there exist strictly greater open convex sets $B$ and $C$ such that $A=B\cap C$.
By the maximality of $A$, both $B$ and $C$ contain $p$, hence $A$ also contains $p$.
This is a contradiction.
We have proved the lemma.
\end{proof}

\begin{theorem}\label{theo:sat-irred}
A subset $S$ of $\p$ is a saturated multi-convex set if and only if $S$ is the intersection of a nonempty family of irreducible convex sets containing $S$.
\end{theorem}

\begin{proof}
Since an irreducible convex set is open, the sufficiency of the condition is satisfied.
The necessity immediately follows from Lemma~\ref{lem:separation}.
\end{proof}

\subsection{A duality for saturated multi-convex sets}
Recall that there exists a natural one-to-one correspondence between irreducible convex sets in $\p$ and points in $\pstar$ (see Proposition~\ref{prop:projduality}).
Given a saturated multi-convex set $C$ in $\p$, this correspondence allows to transform the family of all irreducible convex sets containing $C$ into a set of points in $\pstar$.

\begin{definition}
We define the following function $\Phi\colon 2^{\p}\rightarrow 2^{\pstar}$,
\begin{equation*}
\Phi(S) = \set{\delta(C)\in\pstar}{\text{$C$ is an irreducible convex set in $\p$ containing $S$}},
\end{equation*}
where $\delta$ denotes the one-to-one correspondence between irreducible convex sets in $\p$ and points in $\p^{\ast}$.
In the same way, we have a function from $2^{\p^{\ast}}$ to $2^{\p}$ and denote it by the same symbol $\Phi$ when there is no danger of confusion.
\end{definition}

The following proposition states that $\Phi$ is a Galois connection between the poset $2^{\p}$ ordered by the inclusion and the poset $2^{\p^{\ast}}$ ordered by the reserve inclusion (see Subsection~\ref{subsect:pre}).

\begin{proposition}
\begin{enumerate}
\item The functions $\Phi\colon 2^{\p}\to 2^{\pstar}$ and $\Phi\colon 2^{\pstar}\to 2^{\p}$ are antitone.
\item The relations $S\subseteq \Phi(T)$ and $T\subseteq\Phi(S)$ are equivalent for all pairs $(S,T)\in 2^{\p}\times 2^{\pstar}$.
\end{enumerate}
\end{proposition}

\begin{proof}
The part (i) is trivial. The part (ii) is proved as follows:
\begin{align*}
S\subseteq \Phi(T) & \Longleftrightarrow \text{$T\subseteq\delta(p)$\quad for all points $p\in S$}\\
                   & \Longleftrightarrow \text{$q\in\delta(p)$\quad for all points $p\in S$ and $q\in T$}\\ 
                   & \Longleftrightarrow \text{$p\in\delta(q)$\quad for all points $q\in T$ and $p\in S$}\\ 
                   & \Longleftrightarrow \text{$S\subseteq\delta(q)$\quad for all points $q\in T$}\\ 
                   & \Longleftrightarrow T\subseteq\Phi(S).
\end{align*}
Note that $\delta(p)$ and $\delta(q)$ denote the irreducible convex sets corresponding to $p$ and $q$, respectively.
\end{proof}

\begin{proposition}\label{prop:infsup}
For a family $\{S_{i}\}_{i\in I}$ of subsets of $\p$, we have the following two properties:
\begin{enumerate}
\item $\Phi(\bigcup_{i\in I} S_{i}) =\bigcap_{i\in I}\Phi(S_{i})$;
\item $\Phi(\bigcap_{i\in I} S_{i}) \supseteq\bigcup_{i\in I}\Phi(S_{i})$.
\end{enumerate}
\end{proposition}

\begin{proof}
(i) $p\in\Phi(\bigcup_{i\in I} S_{i})\Leftrightarrow\bigcup_{i\in I}S_{i}\subseteq\delta(p)\Leftrightarrow \text{$S_{i}\subseteq\delta(p)$ for all $i\in I$}\Leftrightarrow p\in\bigcap_{i\in I} \Phi(S_{i})$.

(ii) Since $\Phi$ is antitone, we obtain $\Phi(\bigcap_{i\in I} S_{i})\supseteq\Phi(S_{j})$ for all $j\in I$.
Hence we obtain $\Phi(\bigcap_{i\in I} S_{i})\supseteq\bigcup_{j\in I}\Phi(S_{j})$.
\end{proof}
In the part (ii) above, both sides of the relation do not coincide in general, however if the family of convex sets consists of the co-components of some multi-convex set, then both sides of it coincide.

\begin{proposition}\label{prop:infsup2}
For a multi-convex set $C$ in $\p$, we have 
\[
\Phi(\bigcap\cocomp{C}) =\bigcup\set{\Phi(N)}{N\in\cocomp{C}}.
\]
\end{proposition}

\begin{proof}
By Proposition\ref{prop:infsup}(ii), it suffices to prove $\Phi(\bigcap\cocomp{C})\subseteq\bigcup\set{\Phi(N)}{N\in\cocomp{C}}$.
Let $p$ be a point in $\Phi(\bigcap\cocomp{C})$.
Recall that $\delta(p)$ is an irreducible convex set in $\p$ containing $\bigcap\cocomp{C}$.
Since $C$ is the intersection of its co-components, we obtain $C\subseteq\delta(p)$.
The convex hull $\relhull{C}{\delta(p)}$ of $C$ relative to $\delta(p)$ is a co-component of $C$, which implies that the point $p$ belongs to $\Phi(\relhull{C}{\delta(p)})$.
Hence the point $p$ belongs to $\bigcup\set{\Phi(N)}{N\in\cocomp{C}}$.
\end{proof}

\begin{notation}
Let $\sat{\cdot}$ be an abbreviation for $\Phi\circ\Phi$.
We call the set $\sat{S}$ the {\em saturation} of a subset $S$ of $\p$ (or $\pstar$).
\end{notation}

Since $\Phi$ establishes a Galois connection, we immediately obtain the following result (see Theorem~\ref{theo:galois}).
\begin{proposition}\label{prop:linhull}
For subsets $S$ and $T$ of $\p$, we have the following four properties:
\begin{enumerate}
\item $S\subseteq \sat{S}$;
\item $S\subseteq T$ implies $\sat{S}\subseteq \sat{T}$;
\item $\sat{\sat{S}}=\sat{S}$;
\item $\sat{\Phi(S)}=\Phi(S)$.
\end{enumerate}
\end{proposition}

\begin{proposition}\label{prop:sat-linhull}
For a subset $S$ of $\p$, the following conditions are equivalent:
\begin{enumerate}
\item $S$ is a saturated multi-convex set;
\item $S$ is a proper subset of $\p$ and $\sat{S}=S$.
\end{enumerate}
\end{proposition}

\begin{proof}
We first prove $\Phi\circ\Phi(S)=\bigcap\set{\delta(q)}{q\in\Phi(S)}$:
\begin{align*}
p\in \Phi\circ\Phi(S) & \Longleftrightarrow \Phi(S)\subseteq\delta(p)\\
                      & \Longleftrightarrow \text{$q\in\delta(p)$\quad for all points $q\in\Phi(S)$}\\
                      & \Longleftrightarrow \text{$p\in\delta(q)$\quad for all points $q\in\Phi(S)$}\\
                      & \Longleftrightarrow p\in\bigcap\set{\delta(q)}{q\in\Phi(S)}.
\end{align*}
Note that $\delta(p)$ and $\delta(q)$ denote the irreducible convex sets corresponding to $p$ and $q$, respectively.

Assume the condition (i). Then $S$ is clearly a proper subset of $\p$.
By Theorem~\ref{theo:sat-irred}, the set $S$ is the intersection of all irreducible convex sets containing $S$.
Recall that $\Phi(S)$ is the set of points in $\pstar$ corresponding to the family of all irreducible convex sets containing $S$ in $\p$.
Thus we derive $S=\bigcap\set{\delta(q)}{q\in\Phi(S)}=\Phi\circ\Phi(S)=\sat{S}$.

Conversely assume the condition (ii).
Then we obtain $S=\sat{S}=\Phi\circ\Phi(S)=\bigcap\set{\delta(q)}{q\in\Phi(S)}$.
Hence $S$ is the intersection of a nonempty family of irreducible convex sets.
By Theorem~\ref{theo:sat-irred}, the set $S$ is a saturated multi-convex set.
\end{proof}

Note that we have the corresponding results in $\pstar$ of the preceding four propositions in $\p$, although we omit the statements.

\begin{theorem}
The function $\Phi$ is an order anti-isomorphism between nonempty saturated multi-convex sets in $\p$ and those in $\p^{\ast}$.
\end{theorem}

\begin{proof}
Let $S$ be a nonempty saturated multi-convex set in $\p$.
We prove that $\Phi(S)$ is a nonempty saturated multi-convex set in $\pstar$.
By Proposition~\ref{prop:linhull}(iv), we have $\sat{\Phi(S)}=\Phi(S)$.
Clearly $\Phi(S)$ is a proper subset of $\pstar$.
By the corresponding result in $\pstar$ of Proposition~\ref{prop:sat-linhull}, we deduce that $\Phi(S)$ is a saturated multi-convex set in $\pstar$.
Clearly $\Phi(S)$ is nonempty.
Similarly we can prove that $\Phi\colon 2^{\pstar}\to 2^{\p}$ sends a nonempty saturated multi-convex set in $\pstar$ to a nonempty saturated multi-convex set in $\p$.

By Proposition~\ref{prop:sat-linhull}, we have $\Phi\circ\Phi(S)=S$ for all saturated multi-convex sets $S$ in $\p$ (or $\pstar$).
This implies that $\Phi$ is a bijection between nonempty saturated multi-convex sets in $\p$ and those in $\pstar$.
\end{proof}

\begin{remark}\label{rem:linconv}
In a complex projective space, the notion of linearly convex sets has been studied in \cite{martin:convex} and \cite{znam:convex}.
Let us introduce the notion of linearly convex sets in a real projective space $\p$.
We first define the following function $f$ from $2^{\p}$ to $2^{\pstar}$:
\begin{equation*}
S\mapsto f(S)=\set{w^{\ast}\in\pstar}{\text{$w$ is a hyperplane in $\p$ avoiding $S$}},
\end{equation*}
where $w^{\ast}$ denotes the point in $\pstar$ corresponding to the hyperplane $w$ in $\p$ through the projective duality.
In the same way, we have a function from $2^{\pstar}$ to $2^{\p}$ and denote it by the same symbol $f$.
Then a subset $S$ of $\p$ is {\em linearly convex} if $f\circ f(S)=S$.
By Theorem~\ref{theo:irred}, the function $f$ turns out to coincide with $\Phi$, and by Proposition~\ref{prop:sat-linhull}, a subset of $\p$ is a saturated multi-convex set if and only if it is a linearly convex, proper subset of $\p$.
\end{remark}

\subsection{Components versus co-components}\label{subsect:comp-cocomp}
As expected, the function $\Phi$ preserves the convexity, and furthermore $\Phi$ interchanges the openness and the closedness of convex sets.

\begin{proposition}\label{prop:convtoconv}
For a nonempty subset $S$ of $\p$, we have the following three properties:
\begin{enumerate}
\item if $S$ is a convex set, then $\Phi(S)$ is a convex set;
\item if $S$ is an open convex set, then $\Phi(S)$ is a closed convex set;
\item if $S$ is a closed convex set, then $\Phi(S)$ is an open convex set.
\end{enumerate}
\end{proposition}

\begin{proof}
Let $\p$ be the $n$-dimensional projective space associated with an $(n+1)$-dimensional vector space $V$, and let $\pstar$ be the dual projective space of $\p$.
Let $S$ be a convex set in $\p$.
We can take an affine hyperplane $A$ in $V$ such that $A$ avoids the origin of $V$ and $S\subseteq\pi(A)$, where $\pi$ denotes the projection of $V\setminus\{0\}$ onto $\p$ (see Subsection~\ref{subsect:pre}).
Then $S$ is identified with a convex set in the affine space $A$.

(i) For distinct two points $p,q\in\Phi(S)$, 
let $K$ and $L$ be $1$-dimensional linear subspaces of $V^{\ast}$ such that $\pi(K)=p$ and $\pi(L)=q$.
Note that the orthogonal complements $K^{\bot}$ and $L^{\bot}$ are distinct $n$-dimensional linear subspaces of $V$.
The intersection of $K^{\bot}$ and $L^{\bot}$ is an $(n-1)$-dimensional linear subspace of $V$, hence the set of all $n$-dimensional linear subspaces of $V$ containing $K^{\bot}\cap L^{\bot}$ forms a projective line in $\pstar$.
Among those $n$-dimensional linear subspaces of $V$, the set of all $n$-dimensional linear subspaces of $V$ meeting $S$ is mapped to an interval in $\pstar$.
Since both of $p$ and $q$ lie in the complement of the interval, one of the two line segments joining $p$ and $q$ are contained in $\Phi(S)$, but the other is not.
We have proved the property (i).

(ii) Let $\mathcal{T}(S)$ be the set of all orthogonal complements $W^{\bot}$ such that $W$ is an $n$-dimensional linear subspace of $V$ meeting $S$.
Then $\mathcal{T}(S)$ is a subset of $\pstar$.
Since $S$ is open in the affine space $A$, we can derive that $\mathcal{T}(S)$ is open.
The set $\mathcal{T}(S)$ is the complement of $\Phi(S)$, hence $\Phi(S)$ is closed. 
By the property (i), the set $\Phi(S)$ is convex.

(iii) Note that $\Phi(S)$ coincides with the set of all orthogonal complements $W^{\bot}$ such that $W$ is an $n$-dimensional linear subspace of $V$ avoiding $S$.
Since $S$ is compact in the affine space $A$, we can derive that $\Phi(S)$ is open.
By the property (i), the set $\Phi(S)$ is convex.
\end{proof}

We call a convex set to be {\em saturated} if it is saturated as a multi-convex set.

\begin{corollary}\label{cor:duality3}
The function $\Phi$ is an order anti-isomorphism between nonempty saturated convex sets in $\p$ and those in $\p^{\ast}$ which interchanges the openness and the closedness of convex sets.
\end{corollary}

Moreover, $\Phi$ interchanges the openness and the closedness of multi-convex sets.

\begin{proposition}\label{prop:fliptop}
The function $\Phi$ sends a nonempty open multi-convex set to a closed multi-convex set, and 
a nonempty closed multi-convex set to an open multi-convex set.
\end{proposition}

\begin{proof}
Let $C$ be a nonempty open multi-convex set in $\p$.
We obtain $\Phi(C)=\Phi(\bigcup\comp{C})=\bigcap\set{\Phi(N)}{N\in\comp{C}}$ from Proposition~\ref{prop:infsup}(i).
Since every component $N$ of $C$ is an open convex set, $\Phi(N)$ is a closed convex set.
Thus $\Phi(C)$ is a closed multi-convex set.

Let $C$ be a nonempty closed multi-convex set in $\p$.
We obtain $\Phi(C)=\Phi(\bigcap\cocomp{C})=\bigcup\set{\Phi(N)}{N\in\cocomp{C}}$ from Proposition~\ref{prop:infsup2}.
Then we can similarly prove that $\Phi(C)$ is an open multi-convex set in $\pstar$.
\end{proof}

For a saturated multi-convex set $C$ in $\p$, let us consider the following two statements:
\begin{enumerate}
\item[(1)] every component of $C$ is saturated;
\item[(2)] every co-component of $C$ is saturated.
\end{enumerate}
The statement (1) is clearly true, however the statement (2) is not in general:
for example, the multi-convex set illustrated in Figure~\ref{fig:antiexam-b} is saturated, while its co-components illustrated in Figure~\ref{fig:antiexam-c} and~\ref{fig:antiexam-d} are not saturated.
Thus the function $\Phi$ need not establish a bijection between $\cocomp{C}$ and $\comp{\Phi(C)}$ (and similarly between $\comp{C}$ and $\cocomp{\Phi(C)}$).
However, if $C$ is open or closed, then every co-component of $C$ is open or closed, hence every co-component of $C$ is saturated.
Thus we have:

\begin{figure}[tb]
  \begin{minipage}{0.32\hsize}
    \begin{center}
      \subfigure[]{\includegraphics[height=3.5cm]{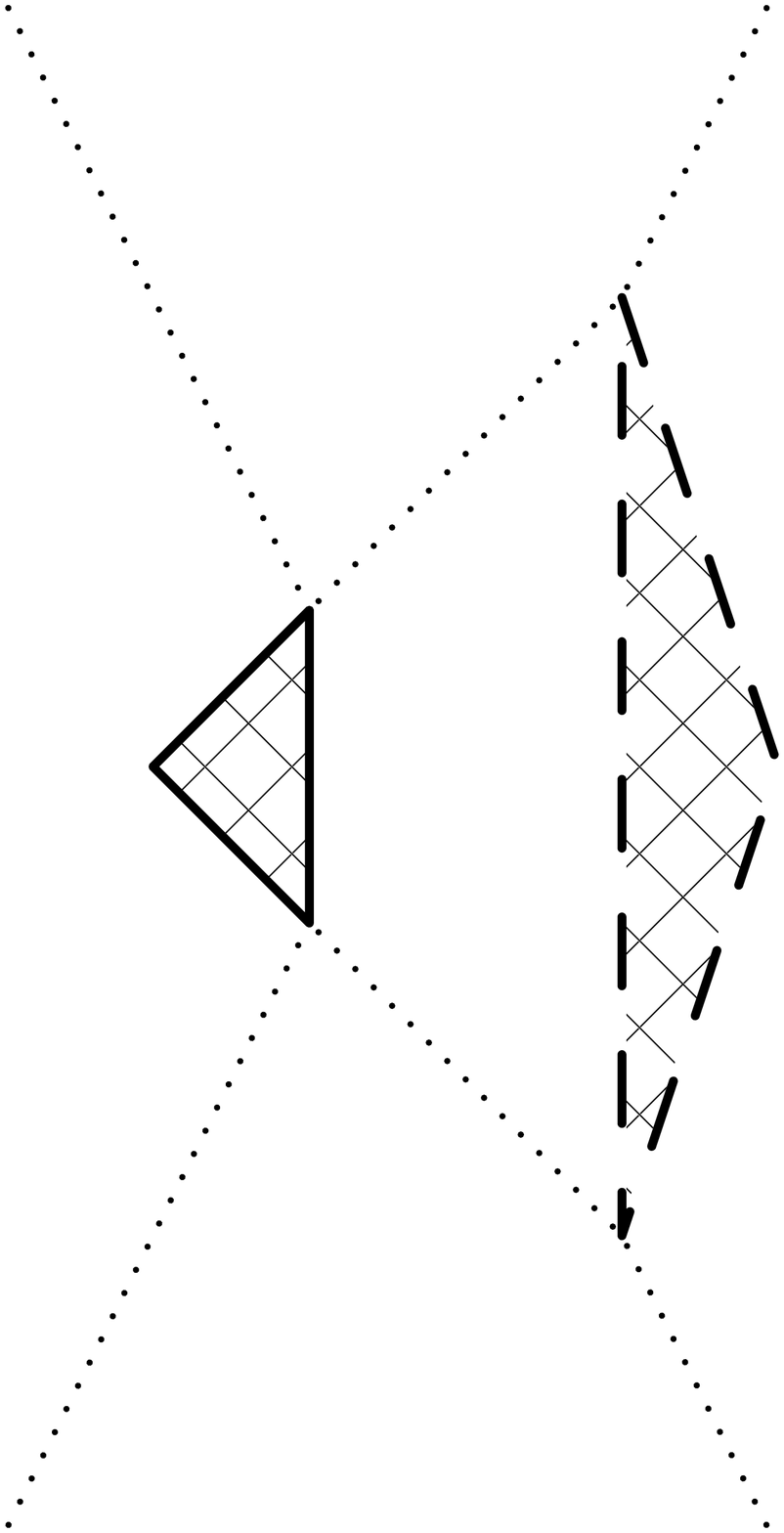}\label{fig:antiexam-b}}
    \end{center}
  \end{minipage}
  \begin{minipage}{0.32\hsize}
    \begin{center}
      \subfigure[]{\includegraphics[height=3.5cm]{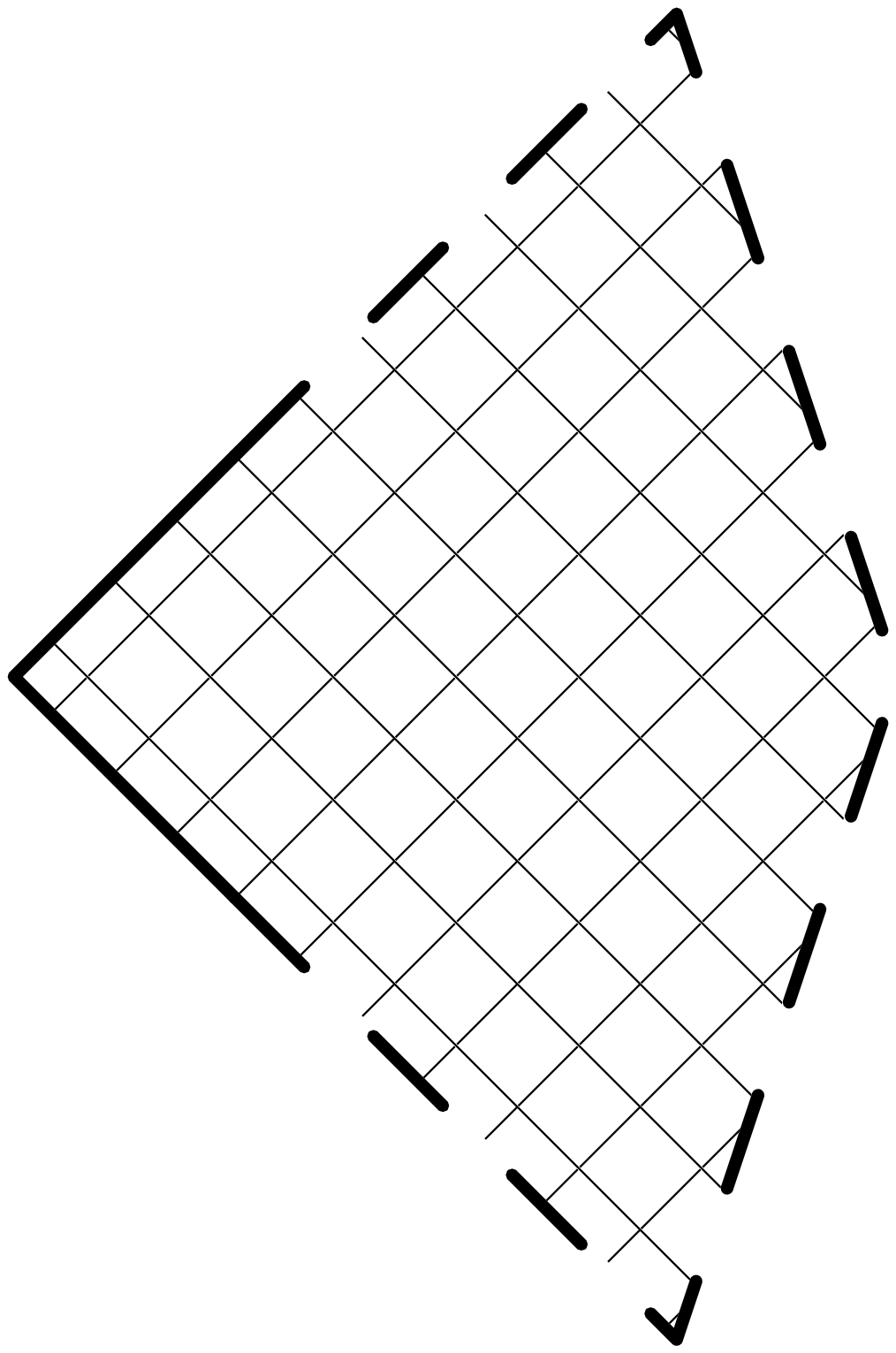}\label{fig:antiexam-c}}
    \end{center}
  \end{minipage}
  \begin{minipage}{0.32\hsize}
    \begin{center}
      \subfigure[]{\includegraphics[height=3.5cm]{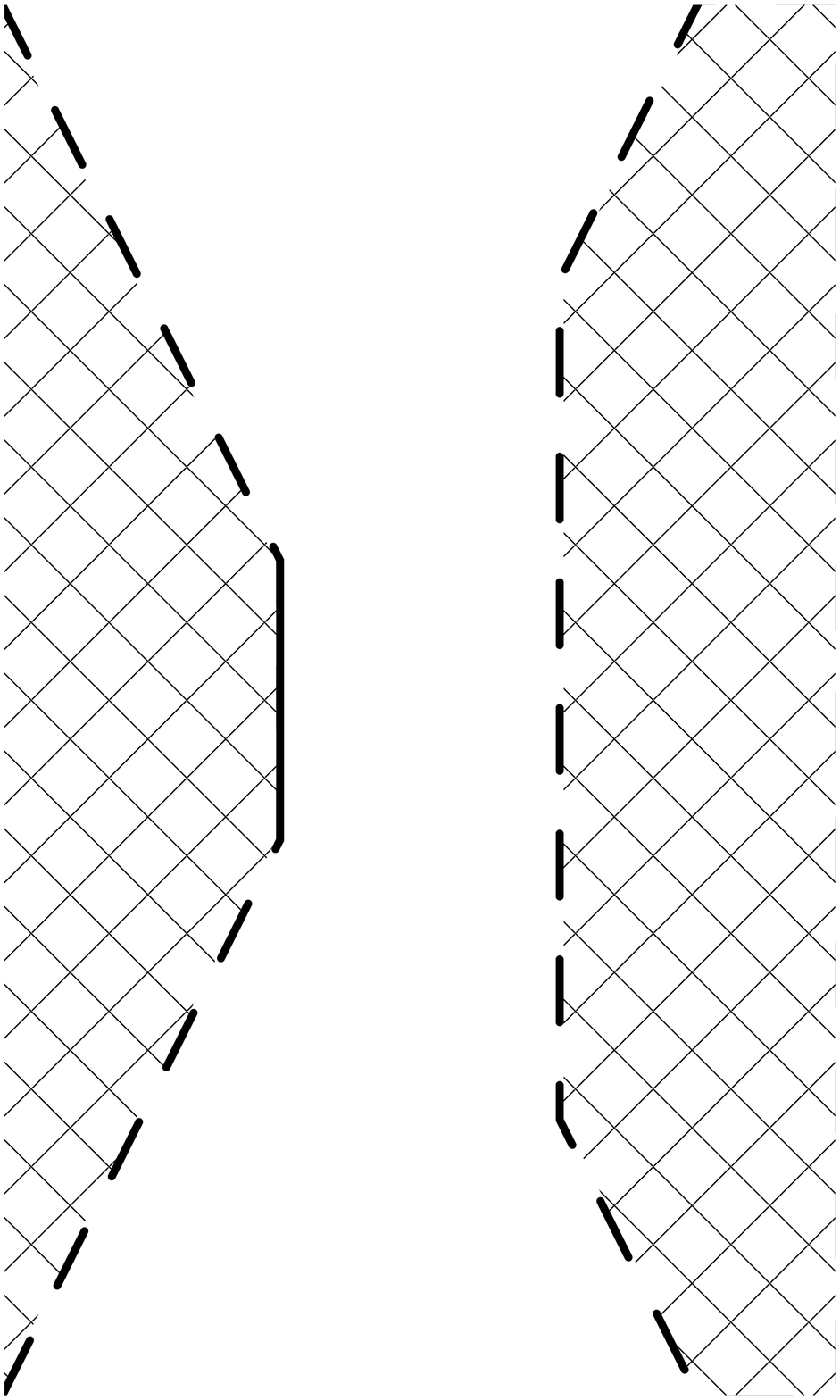}\label{fig:antiexam-d}}
    \end{center}
  \end{minipage}
  \caption{The multi-convex set illustrated in Figure (a) is saturated, while its co-components illustrated in Figure (b) and (c) are not saturated.}\label{fig:antiexam}
\end{figure}

\begin{proposition}
Let $C$ be a nonempty open or closed multi-convex set in $\p$.
Then $\Phi$ is a bijection between $\comp{C}$ and $\cocomp{\Phi(C)}$, and between $\cocomp{C}$ and $\comp{\Phi(C)}$.
\end{proposition}

\begin{proof}
Clearly every component of $C$ is saturated.
By Corollary~\ref{cor:duality3}, the function $\Phi$ is a bijection between components of $C$ and minimal saturated convex sets containing $\Phi(C)$.
By Proposition~\ref{prop:fliptop}, the set $\Phi(C)$ is an open or closed multi-convex set in $\pstar$.
Thus every co-component of $\Phi(C)$ is saturated.
Hence minimal saturated convex sets containing $\Phi(C)$ are the co-components of $\Phi(C)$.
Thus $\Phi$ is a bijection between $\comp{C}$ and $\cocomp{\Phi(C)}$.
Similarly we can prove that $\Phi$ is a bijection between $\cocomp{C}$ and $\comp{\Phi(C)}$.
\end{proof}

For a nonempty saturated multi-convex set $C$, even if $\Phi$ is not a bijection between $\comp{C}$ and $\cocomp{\Phi(C)}$, 
the pair of the following functions derived from $\Phi$ establishes a bijection between them.

\begin{theorem}
Let $C$ be a nonempty saturated multi-convex set in $\p$.
Then the pair of the following functions is a bijection between $\comp{C}$ and $\cocomp{\Phi(C)}$:
\begin{align*}
&\cdot^{\triangleright}\colon \comp{C}\to\cocomp{\Phi(C)},\ N\mapsto\relhull{\Phi(C)}{\Phi(N)};\\
&\cdot^{\triangleleft}\colon \cocomp{\Phi(C)}\to\comp{C},\ M\mapsto\Phi(\sat{M}).
\end{align*}
\end{theorem}

\begin{proof}
We first prove $(N^{\triangleright})^{\triangleleft}=N$ for all components $N$ of $C$.
By Corollary~\ref{cor:duality3}, the set $\Phi(N)$ is a minimal saturated convex set containing $\Phi(C)$.
Clearly it is a saturation of the convex hull $\relhull{\Phi(C)}{\Phi(N)}$ of $\Phi(C)$ relative to $\Phi(N)$, and we obtain
\[
\Phi(N)=\sat{\relhull{\Phi(C)}{\Phi(N)}}=\sat{N^{\triangleright}}.
\]
Since $N$ is saturated, we derive $N=\Phi\circ\Phi(N)=\Phi(\sat{N^{\triangleright}})=(N^{\triangleright})^{\triangleleft}$.

All that remains is to prove $(M^{\triangleleft})^{\triangleright}=M$ for all co-components $M$ of $\Phi(C)$.
We have $\Phi\circ\Phi(\sat{M})=\sat{M}$.
Thus we derive $(M^{\triangleleft})^{\triangleright} = \relhull{\Phi(C)}{\sat{M}}$.
Both $\relhull{\Phi(C)}{\sat{M}}$ and $M$ are co-components of $\Phi(C)$ and are contained in $\sat{M}$.
This implies $\relhull{\Phi(C)}{\sat{M}}=M$. Hence we have $(M^{\triangleleft})^{\triangleright}=M$.
\end{proof}

In the same way, we have a bijection between between $\cocomp{C}$ and $\comp{\Phi(C)}$.

\section*{Acknowledgements}
The author is grateful to Professor Hideki Tsuiki and Professor Klaus Keimel for their encouragement, advices, and comments during the preparation of this paper.
The author also thanks Professor Masahiro Hachimori for providing some important references.

\end{document}